\documentclass[12pt]{article}
\usepackage{amsmath,amssymb,amsbsy,amsfonts,amsthm,latexsym,
                    amsopn,amstext,amsxtra,euscript,amscd,mathrsfs}

\newtheorem{lem}{Lemma}
\newtheorem{lemma}[lem]{Lemma}

\newtheorem{thm}{Theorem}
\newtheorem{theorem}[thm]{Theorem}
\newtheorem{cor}{Corollary}

\def\mand{\qquad \mbox{and} \qquad}

\def\\{\cr}
\def\({\left(}
\def\){\right)}
\def\[{\left[}
\def\]{\right]}
\def\<{\langle}
\def\>{\rangle}

\def\cA{\mathcal A}

\def\fH{\mathfrak H}

\def\fK{\mathfrak K}

\def\fM{\mathfrak M}
\def\cN{\mathcal N}
\def\fN{\mathfrak N}

\def\cP{\mathcal P}
\def\cQ{\mathcal Q}

\def\cS{\mathcal S}

\def\cV{\mathcal V}
\def\cW{\mathcal W}

\def\eps{\varepsilon}

\def\Z{\mathbb{Z}}

\def\r{{\mathscr P}}
\def\paq{{\mathscr P}_{\!a,q}}
\def\vaq{{\mathscr V}_{\!a,q}}
\def\waq{{\mathscr W}_{\!a,q}}

\begin{document}

\title{Integers without divisors\break
from a fixed arithmetic progression}

\author{
{\sc William D.~Banks} \\
{Department of Mathematics, University of Missouri} \\
{Columbia, MO 65211 USA} \\
{\tt bbanks@math.missouri.edu} \\
\and
{\sc John B.~Friedlander} \\
{Department of Mathematics, University of Toronto} \\
{Toronto, Ontario M5S 3G3, Canada} \\
{\tt frdlndr@math.toronto.edu} \\
\and
{\sc Florian~Luca} \\
{Instituto de Matem{\'a}ticas}\\
{ Universidad Nacional Aut\'onoma de M{\'e}xico} \\
{C.P. 58089, Morelia, Michoac{\'a}n, M{\'e}xico} \\
{\tt fluca@matmor.unam.mx}}

\date{}

\pagenumbering{arabic}

\maketitle

\begin{abstract}
Let $a$ be an integer and $q$ a prime number. In this paper we
find an asymptotic formula for the number of positive integers
$n\le x$ with the property that no divisor $d> 1$ of $n$ lies in
the arithmetic progression $a$ modulo $q$.
\end{abstract}



\section{Introduction}

We consider the frequency of natural numbers which do not have any
divisor from a given arithmetic progression. More precisely, for
integers $0\le a<m$ and a real number $x\ge 1$, we define:
$$
\cN(m,a)=\{n\ge 1~:~d\not\equiv a\pmod m\text{~for all~}d\mid n,~d
>1\}
$$
and denote by $\cN(x;m,a)$ the number of positive integers $n\le
x$ in $\cN(m,a)$.

Our goal is to determine an explicit asymptotic formula for
$\cN(x;m,a)$. We exclude the divisor $d=1$ in the above definition
since including it would make the result trivial for one residue
class while not affecting the result for any of the others. To
avoid increasing the technical complications we give detailed
consideration to the special case that $m=q$ is a prime number. In
the final section we give some remarks about the case of general
modulus and about the still more complicated problem of counting
those integers whose divisors avoid a subset of the residue
classes.

When $a=0$, it is clear that $n\in\cN(q,0)$ if and only if $q$
does not divide~$n$, and in this case it follows that
$$
\cN(x;q,0)=(1-q^{-1})\,x+O(q).
$$
Thus,  we can assume that $a\ge 1$ in what follows.

If $a=1$ and $q=2$, it is also clear that $n$ is in $\cN(2,1)$ if
and only if $n$ is a power of two, and therefore,
$$
\cN(x;2,1)=\frac{\log x}{\log 2}+O(1).
$$
Hence, we can further assume that $q\ge 3$ throughout the sequel.
The case $a=1$ is essentially different from (and quite a bit
easier than) the others. The result obtained is the following.
\begin{theorem}
\label{thm:easy} For every fixed odd prime $q$ we have
$$
\cN(x;q,1)=(1+o(1))\,\frac{\varphi(q-1)\,q}{(q-1)^{q-1}\,(q-3)!}
\,\frac{x(\log\log x)^{q-3}}{\log x}\ ,
$$
where $\varphi$ is the Euler function.
\end{theorem}

In view of Theorem~\ref{thm:easy}, which is proved in
Section~\ref{sec:thm1} below, it remains only to consider the case
that $1<a<q$.  In order to state this result we introduce three
constants $\paq$, $\vaq$, and $\waq$, as follows. First, let
$$
q-1=\prod_{j=1}^{k} p_j^{\alpha_j}\qquad {\text{\rm and}}\qquad
{\rm ord\/}_q(a)=\prod_{j=1}^s p_j^{\beta_j}
$$
be the prime factorizations of $q-1$ and ${\rm ord\/}_q(a)$ (the
multiplicative order of $a$ modulo $q$), respectively. Here,
$p_1,\ldots,p_k$ are distinct primes, $s\le k$, and the integers
$\alpha_j$ and $\beta_j$ are positive. Using these data, we define
\begin{equation}
\label{eq:paq_defn} \paq=\min_{1\le j\le
s}\big\{p_j^{\alpha_j-\beta_j+1}\big\}.
\end{equation}
Next, recalling that every subgroup of a cyclic group is
determined uniquely by its cardinality, let $H(a)$ be the unique
subgroup in $(\Z/q\Z)^*$ of cardinality $|H(a)|=(q-1)/\paq$, and
put
\begin{equation}
\label{eq:vaq_defn} \vaq=\lim_{y\to\infty}\Biggl\{(\log
y)^{-1/\paq}\!\!\!\!\!\!
   \prod_{\substack{p\le y\\p\!\!\!\pmod q\in
H(a)}}\!\!\!\!\!\!\(1-\frac{1}{p}\)^{-1}\Biggl\}.
\end{equation}
Thanks to the work of Williams~\cite{Williams}, one knows that the
limit exists and $0<\vaq<\infty$. Finally, suppose that $\paq$ is
the prime power $p^r$, and put
\begin{equation}
\label{eq:waq_defn}
\waq=\frac{e^{-\gamma/\paq}(1-q^{-1})^{-1+1/\paq}}
{\Gamma(1/\paq)\paq^{\paq-2}\,(\paq-2)!}~p^{r-1}\sum_{1\le j\le
p/2}\binom{p^r-2}{p^{r-1}j-1},
\end{equation}
where $\gamma$ is the \emph{Euler-Mascheroni constant}, and $
\Gamma(s)$ is the usual gamma function.

\begin{theorem}
\label{thm:hard} For every fixed odd prime $q$ and integer $a$
with $1<a<q$, we have
$$
\cN(x;q,a)=(1+o(1))\,\vaq\waq\,\frac{x(\log\log x)^{\paq-2}}{(\log
x)^{1-1/\paq}}\ .
$$
\end{theorem}

We deal throughout with a fixed arithmetic progression and do not
consider the question of uniformity of the estimations in the
modulus $q$, although it is clear from the methods employed that
some (probably not very large) range of uniformity could be
obtained.

The question of counting the number of integers up to $x$ with no
{\it prime} divisor in a given residue class is more familiar and
has a simpler answer; see for example the theorem of Wirsing given
below in Lemma~\ref{lem:wirs}. Our proofs use this result and
similar analytic methods but are complicated by other
considerations which are mostly of a combinatorial nature and with
a bit of group theory.

As we shall see in Lemma~\ref{lem:G/H(a)}, the group $H(a)$ is the
subgroup of $(\Z/q\Z)^*$ having the largest order amongst those
which do not contain (the class of) $a$, and this suggests its
relevance to our problem. The fact that this subgroup is not
unique in general, when the group is not cyclic, is the main thing
which complicates the case of arbitrary modulus. These facts also
lead, in our case of prime modulus, to the following easy
corollaries.
\begin{cor}\label{cor:equigroup}
In case $H(a)=H(b)$ we have
$$
\cN(x;q,a)\sim \cN(x;q,b)\ .
$$
\end{cor}
Special cases of this give the following two results.
\begin{cor}\label{cor:conj}
If $\overline a$ satisfies $a \overline a \equiv 1 \pmod q$ then
$$
\cN(x;q,\overline a)\sim \cN(x;q,a)\ .
$$
\end{cor}
\begin{cor}\label{cor:nonres}
If $a$ and $b$ are both quadratic non-residues modulo $q$ then
$$
\cN(x;q,a)\sim \cN(x;q,b)\ .
$$
\end{cor}
\begin{proof}
In this case $H(a)$ and $H(b)$ are each the subgroup of quadratic
residues.
\end{proof}
Finally we have
\begin{cor}\label{cor:resnonres}
If $a$ is a quadratic residue modulo $q$ and $b$ is a quadratic
non-residue then
$$
\cN(x;q,a)=o(\cN(x;q,b))\ .
$$
\end{cor}
\begin{proof}
In this case, either $a=1$ and the result follows on comparing the
estimates of the two theorems or, if $a>1$, then $H(a)$ is a
subgroup of index greater than two and the result follows from the
second theorem.
\end{proof}

\medskip

Although there seem to be no earlier results that consider the
above asymptotic formulae in this rather basic question, there is
a long history of work on closely related problems. Erd\H
os~\cite{E} showed that, if $m\le (\log x)^{\log 2-\delta}$ where
$\delta>0$ is fixed, then almost all positive integers $n\le x$
have a divisor $d$ in each one of the residue classes $a\pmod m$,
with $\gcd(a,m)=1$. The value $\log 2$ is optimal. Indeed, if $n$
satisfies the above condition then $\tau(n)\ge \phi(m)$ and, since
$\tau(n)=(\log x)^{\log 2+o(1)} $ holds for almost all $n\le x$,
we find that $m\le (\log x)^{\log 2+o(1)}$. Since the appearance
of~\cite{E}, the distribution of integers having a divisor in a
specific residue class has been studied by several authors. For
example, in answer to a question of Erd\H os from~\cite{E1},
Hall~\cite{Ha} showed that, for any $\varepsilon>0$ and natural
number $N$, there exists $\eta_N$ with $\eta_N\to 0$ when
$N\to\infty$ such that, if $m\ge (\log N)^{\log 2}
2^{\sqrt{(2+\varepsilon)\log\log N \log\log\log\log N}}$ then the
number of positive integers $n\le x$ having a divisor $d$ in the
interval $m\le d\le N$ with $d\equiv 1\pmod m$ is $<\eta_N x$
provided $N\le x$. Extending prior results of Hall~\cite{Ha} and
Erd\H os and Tenenbaum~\cite{ET}, de la Bret\'eche~\cite{B} proved
that, if $N$ is any positive integer and $z_N$ is defined
implicitly by the relation $m=(\log N)^{\log 2}
2^{-z_N{\sqrt{\log\log N}}}$, then there exists $\eta_N\to 0$ when
$N\to\infty$ such that, for any $a$ coprime to $m$, we have
$$
\Big|\{n\le x~:~d\equiv a\pmod m~{\text{\rm for~some}}~d\mid n,~
m\le d\le N\}\Big|=\Phi(z_N)+O(\eta_N x)
$$
for all $3\le N\le x$, where
$$
\Phi(z)=\frac{1}{{\sqrt{2\pi }}}\int_{-\infty}^{z}e^{-\tau^2/2}
d\tau,
$$
which in turn answered a conjecture of Hall from \cite{Ha1}.

\medskip

Throughout the paper, $x$ denotes a large positive real number. We
use the Landau symbols $O$ and $o$, as well as the Vinogradov
symbol $\ll$, with their usual meanings. As we do not consider the
question of obtaining bounds which are uniform in the modulus of
the arithmetic progression we allow the implied constants in many
places to depend on various parameters, such as the modulus,
without explicit mention.

For a positive integer $\ell$, we write $\log_\ell x$ for the
function defined inductively by $\log_1 x=\max\{\log x,1\}$ and
$\log_\ell x=\log_1(\log_{\ell-1} x)$ for $\ell\ge 2$, where
$\log$ denotes the natural logarithm function. In the case
$\ell=1$, we omit the subscript to simplify the notation; however,
it should be understood that all the logarithms that appear are
$\ge 1$.

We use various other standard notations, including those for basic
arithmetic functions such as Euler's $\varphi$--function. We use
$|A|$ to denote the number of elements in $A$ when $A$ is a finite
group, or set, or multiset. Given a set $S$ of positive integers,
whether finite or infinite, we frequently denote by $S(x)$ the
number of integers $n\le x$ in $S$.

\medskip

{\bf Acknowledgements.} We thank R\'egis de la Bret\'eche for
useful discussions and for supplying us with several references.
Most of the early work on this paper was done during a visit by
the second author to the Mathematical Institute of the UNAM in
Morelia, Mexico, and during a visit of the third author at the
University of Missouri, Columbia. Most of the later work was done
during visits by all three authors to BIRS in Banff and to CRM in
Montreal. The hospitality and support of all these institutions
are gratefully acknowledged. During the preparation of this paper,
W.~B.\ was supported in part by NSF grant DMS-0070628, J.~F.\ was
supported in part by  NSERC grant A5123 and a Killam Research
Fellowship, and F.~L. was supported in part by Grants SEP-CONACyT
46755, PAPIIT IN104005 and a Guggenheim Fellowship.

\section{Preliminary results}
\subsection{Combinatorial results}
\label{sec:comb}

Recall that a \emph{multiset} is a list $\langle
a_1,\ldots,a_k\rangle$ of elements from a set in which the same
element can occur more than once, but the order is unimportant.
For example, $\langle 1,1,2,3\rangle$ and $\langle 3,1,2,1\rangle$
are the same multisets in $\Z$, whereas $\langle 1,1,2,3\rangle$
and $\langle 1,2,3\rangle$ are different.

Let $G$ be an arbitrary finite abelian group, written additively.
If $G=\{0\}$, put $\kappa(G)=0$; otherwise, let $\kappa(G)$ be the
largest integer $k$ for which there exists a multiset $\langle
a_1,\ldots,a_k\rangle$ of elements of $G$ with the property:
\begin{equation}
\label{eq:zerosum} \sum_{j\in S} a_j\ne 0 \quad\text{for every
nonempty subset~}S\subseteq \{1,2,\ldots,k\}.
\end{equation}
Since $|G|<\infty$, it is easy to see that $\kappa(G)<\infty$. In
the special case that $G=\Z/m\Z$, we have the following result:

\begin{lemma}
\label{lem:comb1} Let $G=\Z/m\Z$, where $m\ge 1$. Then
$\kappa(G)=m-1$. Moreover, if $m\ge 2$, then the multiset $\langle
a_1,\ldots,a_{m-1}\rangle$ has the property~\eqref{eq:zerosum} if
and only if $a_1=\cdots=a_{m-1}=a$ for some $a\in G$ that is
coprime to $m$.
\end{lemma}

\begin{proof}
We can assume that $m\ge 2$ since the result is trivial for $m=1$.

Suppose that $\kappa(G)\ge m$.  Then, for some $k\ge m$, there
exists a multiset $\langle a_1,\ldots,a_k\rangle$ in $G$ with the
property~\eqref{eq:zerosum}. Since the elements $b_j=\sum_{i=1}^j
a_i$, $j=1,\ldots,k$, are all nonzero, and $G$ has only $m-1$
nonzero elements, two of the elements $b_j$ must be equal by the
pigeonhole principle; that is, $b_{j_1}=b_{j_2}$ for some
$j_1<j_2$. But this implies that $\sum_{j_1<i\le j_2} a_i=0$,
which contradicts~\eqref{eq:zerosum}. Therefore,
    $\kappa(G)\le m-1$.

Next, suppose that $a_1=\cdots=a_{m-1}=a$ for some $a\in G$ that
is coprime to $m$. Then, for every nonempty subset $S\subseteq
\{1,2,\ldots,m-1\}$, one has $\sum_{j\in S}a_j=a|S|$.  Since
$\gcd(a,m)=1$ and $m\nmid |S|$, it cannot be true that $a|S|\equiv
0\pmod m$; therefore, the multiset $\langle
a_1,\ldots,a_{m-1}\rangle$ has the property~\eqref{eq:zerosum}
which shows that $\kappa(G)\ge m-1$.

Finally, suppose that the multiset $\langle
a_1,\ldots,a_{m-1}\rangle$ has the property~\eqref{eq:zerosum}. As
before, let $b_j=\sum_{i=1}^j a_i$, $j=1,\ldots,m-1$. Then the
elements $b_1,\ldots,b_{m-1}$ are distinct and nonzero, and since
$G$ has precisely $m-1$ nonzero elements, it follows that
$\{b_1,\ldots,b_{m-1}\}=G\setminus\{0\}$.
Using~\eqref{eq:zerosum}, we see that $a_2\ne 0$, and $a_2\ne b_j$
for $j=2,\ldots,m-1$; therefore, $a_2=b_1=a_1$. By a similar
argument, it follows that $a_j=a_1$ for $j=2,\ldots,m-1$; in other
words, $a_1=\cdots=a_{m-1}=a$ holds for some $a\in G$. Thus, we
have $b_j=ja$ for $j=1,\ldots,m-1$ and, since $b_j\equiv 1\pmod m$
for some value of $j$, it follows that $\gcd(a,m)=1$.
\end{proof}

Now, let $G$ be a nontrivial finite abelian group, written
additively. If $G=\Z/2\Z$, put $\kappa(G,1)=0$; otherwise, for
every $a\in G\setminus\{0\}$, let $\kappa(G,a)$ be the largest
integer $k$ for which there exists a multiset $\langle
a_1,\ldots,a_k\rangle$ of elements of $G\setminus\{0\}$ with the
property:
\begin{equation}
\label{eq:zerosum1} \sum_{j\in S} a_j\ne a \quad\text{for every
subset~}S\subseteq \{1,2,\ldots,k\}.
\end{equation}
In general, $\kappa(G,a)$ need not be finite (e.g., if $p$ is
prime, $G=\Z/p^2\Z$, and $a=1$, then~\eqref{eq:zerosum1} holds for
the multiset $\langle a_1,\ldots,a_k\rangle$, where
$a_1=\cdots=a_k=p$, for every natural number $k$).  However, we do
have the following finiteness result, which suffices for our
applications:

\begin{lemma}
\label{lem:comb2} Let $p$ be a fixed prime, and let $G_r=\Z/p^r\Z$
for every natural number $r$.  If $a\in G_r\setminus\{0\}$ and
$p^{r-1}\mid a$, then $\kappa(G_r,a)=p^r-2$.
\end{lemma}

\begin{proof}
First, we argue by induction on $r$ that for every $a\in
G_r\setminus\{0\}$ with $p^{r-1}\mid a$ and every multiset
$\langle a_1,\ldots,a_k\rangle$ in $G_r\setminus\{0\}$ with the
property~\eqref{eq:zerosum1}, the following inequality holds:
\begin{equation}
\label{eq:boundineq} \Big |\Big\{g\in G_r~:~g=\sum_{j\in
S}a_j\text{~for some subset~}S\subseteq\{1,2,\ldots,k\}\Big\}\Big
|\ge k+1.
\end{equation}
Since the left side of~\eqref{eq:boundineq} cannot exceed
$\big|G_r\setminus\{a\}\big|=p^r-1$, it follows that
$\kappa(G_r,a)\le p^r-2$.

Suppose first that $r=1$, and put $G=G_1=\Z/p\Z$. Let
     $a\in G\setminus\{0\}$ be fixed,
and suppose that $\langle a_1,\ldots,a_k\rangle$ is a multiset in
$G\setminus\{0\}$ with the property~\eqref{eq:zerosum1}. Let
$b_1,\ldots,b_s$ be the distinct values taken by $a_i$ for
$i=1,\ldots k$, and let $m_1,\ldots,m_s$ be the respective
multiplicities; then $\sum_{j=1}^s m_j=k$. Put
$$
A_j=\{u b_j~:~u=0,1,\ldots,m_j\}\qquad (1\le j\le s).
$$
Since each $m_j<p-1$ (otherwise, $a\in A_j$
and~\eqref{eq:zerosum1} fails), $A_j$ is a subset of $G$ of
cardinality $m_j+1$. Let $\sum_{j=1}^sA_j$ be the set of elements
$g\in G$ of the form $g=\sum_{j=1}^sc_j$, where $c_j\in A_j$ for
$j=1,\ldots,s$. A corollary/generalization of the {\it
Cauchy-Davenport theorem\/} (see for
example~\cite[Theorem~2.3]{Nat}) states that
$$
\Big|\,\sum_{j=1}^sA_j\Big|\ge\min\Big\{p,\sum_{j=1}^s|A_j|-s+1\Big\},
$$
and in our situation,
$$
\sum_{j=1}^s |A_j|-s+1=\sum_{j=1}^s(m_j+1)-s+1=k+1.
$$
Since $\sum_{j=1}^sA_j$ is the set of elements $g\in G$ that can
be written as $\sum_{j\in S}a_j$ for some subset
$S\subseteq\{1,2,\ldots,k\}$, we also have by~\eqref{eq:zerosum1}:
$$
\Big|\,\sum_{j=1}^sA_j\Big|\le \big|G\setminus\{a\}\big|=p-1.
$$
Therefore, $k+1\le p-1$, and we obtain the
inequality~\eqref{eq:boundineq} when $r=1$.

To complete the induction, we show that~\eqref{eq:boundineq} holds
for the integer $r\ge 2$ assuming that the corresponding
inequality is true for $r-1$.

Let $a\in G_r\setminus\{0\}$ with $p^{r-1}\mid a$, and suppose
that $\langle a_1,\ldots,a_k\rangle$ is a multiset in
$G_r\setminus\{0\}$ satisfying~\eqref{eq:zerosum1}. Without loss
of generality, we can assume that
     $a_1,\ldots,a_\ell\in G_r\setminus G'_r$ and
$a_{\ell+1},\ldots,a_k\in G'_r\setminus \{0\}$, where $G'_r$ is
the subgroup of $G_r$ consisting of those elements divisible by
$p$.

Let $\widetilde a_j=a_{\ell+j}/p$ for $j=1,\ldots,k-\ell$, and put
$\widetilde a=a/p$. Then $\widetilde a\in G_{r-1}\setminus\{0\}$
with $p^{r-2}\mid\widetilde a$, and $\langle\widetilde
a_1,\ldots,\widetilde a_{k-\ell}\rangle$ is a multiset in
$G_{r-1}\setminus\{0\}$ that satisfies the analogous statement
of~\eqref{eq:zerosum1} obtained after replacing $a$ by $\widetilde
a$, each $a_j$ by $\widetilde a_j$, and $k$ by $k-\ell$, since the
condition $\sum_{j\in S} \widetilde a_j\ne \widetilde a$ in
$G_{r-1}$ is equivalent to $\sum_{j\in S} a_{\ell+j}\ne a$ in
$G_r$ for every subset $S\subseteq \{1,2,\ldots,k-\ell\}$.
Applying the inductive hypothesis with the element $\widetilde a$
and the multiset $\langle\widetilde a_1,\ldots,\widetilde
a_{k-\ell}\rangle$ in $G_{r-1}$, and considering its implication
for the element $a$ and the multiset $\langle
a_{\ell+1},\ldots,a_k\rangle$ in $G_r$, one sees that if $B$
denotes the set of elements $g\in G_r$ equal to $\sum_{j\in S}a_j$
for some subset $S\subseteq\{\ell+1,\ldots,k\}$, then
$\big|B\big|\ge k-\ell+1$.

Let $b_1,\ldots,b_s\in G_r\setminus G'_r$ be the distinct values
taken by $a_i$ for $i=1,\ldots,\ell$, and let $m_1,\ldots,m_s$ be
the respective multiplicities; then $\sum_{j=1}^s m_j=\ell$. Let
$A_j=\{0,b_j\}$, and put
$$
m_jA_j=\mathop{\underbrace{A_j+\cdots+A_j}}\limits_{\text{$m_j$
copies}}\qquad (1\le j\le s).
$$
Since each $b_j$ is \emph{coprime} to $p$, a theorem of I.~Chowla
(see~\cite[Theorem~2.1]{Nat}) yields the inequality
$$
\Big|B+\sum_{j=1}^s m_jA_j\Big|\ge\min\Big\{p^r,
\big|B\big|+\sum_{j=1}^s m_j|A_j|-\sum_{j=1}^s(m_j-1)-s\Big\}.
$$
Since $\big|A_j\big|=2$ for $j=1,\ldots,s$, we have
$$
\big|B\big|+\sum_{j=1}^s m_j\big|A_j\big|-\sum_{j=1}^s(m_j-1)-s
\ge k-\ell+1+\sum_{j=1}^s m_j=k+1.
$$
As $B+\sum_{j=1}^s m_jA_j$ is the set of elements $g\in G_r$ that
are equal to $\sum_{j\in S}a_j$ for some subset
$S\subseteq\{1,2,\ldots,k\}$, we also have by~\eqref{eq:zerosum1}:
$$
\Big|B+\sum_{j=1}^s
m_jA_j\Big|\le\big|G_r\setminus\{a\}\big|=p^r-1.
$$
Therefore, $k+1\le p^r-1$, and we obtain the
inequality~\eqref{eq:boundineq}, which completes the induction.

As mentioned earlier, the inequality~\eqref{eq:boundineq} implies
that $\kappa(G_r,a)\le p^r-2$ for all $a\in G_r\setminus\{0\}$
with $p^{r-1}\mid a$. On the other hand, the lower bound
$\kappa(G_r,a)\ge p^r-2$ is an immediate consequence of the next
lemma.
\end{proof}

\begin{lemma}
\label{lem:construct1} Suppose that $p$, $r$, and $a$ satisfy the
conditions of Lemma~\ref{lem:comb2}, and put $k=p^r-2$. For every
$b\in G_r$ such that $p\nmid b$, let $n$ be the least nonnegative
integer for which the congruence $n\equiv ab^{-1}-1\pmod {p^r}$
holds, and let $\fM_{p,r,a}(b)=\langle a_1,\ldots,a_k\rangle$ be
the multiset in $G_r\setminus\{0\}$ defined by
$$
a_j=\begin{cases}
\, \, \, b&\quad\text{if $j\le n$}; \\
-b&\quad\text{if $j\ge n+1$}.
\end{cases}
$$
Then $\fM_{p,r,a}(b)$ has the property~\eqref{eq:zerosum1}.
\end{lemma}

\begin{proof}
For every subset $S\subseteq\{1,\ldots,k\}$, we have $\sum_{j\in
S}a_j=mb$ for some integer $m$ in the range $-(k-n)\le m\le n$.
Hence, $m\not\equiv (n+1)\pmod{p^r}$, and therefore
$mb\not\equiv(n+1)b\equiv a\pmod{p^r}$.
\end{proof}

The next lemma shows that the multisets $\fM_{p,r,a}(b)$ defined
in Lemma~\ref{lem:construct1} are the only critical multisets that
arise under the conditions of Lemma~\ref{lem:comb2}.

\begin{lemma}
\label{lem:extremecase} Suppose that $p$, $r$, and $a$ satisfy the
conditions of Lemma~\ref{lem:comb2}, and put $k=p^r-2$. If
$\langle a_1,\ldots,a_k\rangle$ is a multiset in
$G_r\setminus\{0\}$ with the property~\eqref{eq:zerosum1}, then
$\langle a_1,\ldots,a_k\rangle=\fM_{p,r,a}(b)$ for some choice of
$b\in G_r$.
\end{lemma}

\begin{proof}
We proceed by induction on $r$, following the proof of
Lemma~\ref{lem:comb2}.

First, let $r=1$. Suppose there exist integers $b,c$ with
$b\not\equiv \pm c\pmod p$ and indices $i,j$ such that $a_i\equiv
b\pmod p$ and $a_j\equiv c\pmod p$. Reordering the elements
$a_1,\ldots,a_k$ if necessary, we can assume that $i=1$ and $j=2$.
Let $A=\{0,a_1\}+\{0,a_2\}$; clearly, $|A|=4$. Let $A_j=\{0,a_j\}$
for $j=1,\ldots,k$. By the Cauchy-Davenport theorem, we have
$$
p-1\ge \Big|A+\sum_{j=3}^k A_j\Big| \ge \min\Big\{p,\big|A\big|
+\sum_{j=3}^k \big|A_j\big|-(k-1)+1\Big\}=p,
$$
which is impossible. Thus, there exists an integer $b$ such that
$a_j\in \{b,-b\}$ for $j=1,\ldots, k$. After reordering the
elements $a_1,\ldots,a_k$, we can assume that $a_j=b$ if $j\le m$
and $a_j=-b$ if $j\ge m+1$, for some $0\le m\le k$.

Now, let $n$ be the least positive integer for which $n\equiv
ab^{-1}\pmod{p}$ holds.  If $n\le m$, then
$a_1+\cdots+a_n=nb\equiv a\pmod p$, which
contradicts~\eqref{eq:zerosum1}. On the other hand, if $n\ge m+2$,
then $p-n\le p-2-m=k-m$, thus
$a_{m+1}+\cdots+a_{m+p-n}=(p-n)(-b)\equiv a\pmod p$, which again
contradicts~\eqref{eq:zerosum1}.  Therefore, $n=m+1$, and the
result follows for $r=1$.

Now suppose that the result has been proved for all cyclic
$p$-groups of order less than $p^r$; we need to prove it for
$G_r=\Z/p^r\Z$.

To do this, let us first show that $p\nmid a_j$ for each
$j=1,\ldots,k$. Indeed, suppose this is not the case. If $p\mid
a_j$ for \emph{all} $j=1,\ldots,k$, then writing
$\widetilde{a}_j=a_j/p$, we see that the multiset
$\langle\widetilde{a}_1,\ldots, \widetilde{a}_k\rangle$ has the
property~\eqref{eq:zerosum1} with $a$ replaced by
$\widetilde{a}=a/p$. Since the elements $\widetilde{a}_1,\ldots,
\widetilde{a}_k$ can be viewed as elements of the cyclic group
with $p^{r-1}$ elements, the induction hypothesis implies that
$p^r-2=k\le p^{r-1}-2$, which is impossible. This argument shows
that there exists at least one element $a_j$ such that $p\nmid
a_j$. Now, using the notation of the proof of
Lemma~\ref{lem:comb2}, we can assume that $p\mid a_j$ for
$j=1,\ldots,\ell$ and $p\nmid a_j$ for $j=\ell+1,\ldots,k$, where
$1\le\ell<k$. Let $B$ denote the set of elements $g\in G_r$ that
are equal to $\sum_{j\in S}a_j$ for some subset
$S\subseteq\{\ell+1,\ldots,k\}$. As in the proof of
Lemma~\ref{lem:comb2}, we have $\big|B\big|\ge k-\ell+1$. Since
$p\mid a_1$, it follows that
$$
\big|B+\{0,a_1\}\big|=2\big|B\big|>\big|B\big|+1.
$$
Since the last inequality is strict, the argument based on the
Chowla Theorem (see the proof of Lemma~\ref{lem:comb2}) implies
that
$$
\Big|B+\sum_{j=1}^s m_jA_j\Big|>\min\{p^r,k+2\}.
$$
Taking into account that $B+\sum_{j=1}^s m_jA_j$ has at most
$p^r-1$ elements (since this set does not contain $a$), we see
that $k\le p^r-3$, which is impossible. Thus, we have established
our assertion that $p\nmid a_j$ for $j=1,\ldots,k$.

To complete the proof of the lemma for the group $G_r$, we can use
an argument identical to the one given above for the case $r=1$,
except that the Cauchy-Davenport Theorem is now replaced by the
Chowla Theorem, which is applicable since $\gcd(a_j,p)=1$ for
$j=1,\ldots,k$.
\end{proof}

The next lemma provides a complete list of the \emph{distinct}
multisets $\fM_{p,r,a}(b)$ which arise for various choices of $p$
and $r$ in the special case that $a=p^{r-1}$.

\begin{lemma}
\label{lem:listmultisets} Let $p^r$ be a prime power, and let
$a=p^{r-1}$. For integers $\eta$ and $c$ let $\fN_{p,r}(\eta,c)$
be the multiset in $\Z$ defined by
$$
\fN_{p,r}(\eta,c)=\big\langle\mathop{\underbrace{c,c,\ldots,c}}
\limits_{\text{$\eta$
copies}}~,\,\mathop{\underbrace{-c,-c,\ldots,-c}}
\limits_{\text{$p^r-2-\eta$ copies}}\big\rangle .
$$
For an integer $\lambda$ not divisible by $p$ let
$\overline\lambda$ be the least positive integer such that
$\lambda\overline\lambda\equiv 1\pmod p$. Then, there is a one to
one correspondence, given by the congruence modulo $p^r$, between
pairs of multisets $\fM_{p,r,a}(\pm b)$ in $G_r\setminus\{0\}$ and
the family $\fN_{p,r}(\eta,c)$ where, in case $p$ is odd, $\eta$
runs through the integers
$$\eta\in\{p^{r-1}\lambda-1:1\le\lambda \le(p-1)/2\}
$$
and $c$ runs through the integers
$$c\in\{\overline\lambda+p\mu:0\le\mu\le p^{r-1}-1\}\ ,
$$
while, in case $p=2$, we have $\eta = 1$ and, in the range for
$c$, we must replace $p^{r-1}-1$ by $2^{r-2}-1$.
\end{lemma}

\begin{proof}
Let $\fM_{p,r,a}(b)$ be a multiset in $G_r\setminus\{0\}$ of the
type constructed in Lemma~\ref{lem:construct1}.  We claim that
$\fM_{p,r,a}(b)=\fM_{p,r,a}(-b)$. Indeed, let $n$ be the least
nonnegative integer for which the congruence $n\equiv
p^{r-1}b^{-1}-1\pmod{p^r}$ holds. Clearly, $n\ne p^r-1$, hence it
follows that $m=p^r-2-n$ is the least nonnegative integer for
which the congruence $m\equiv p^{r-1}(-b)^{-1}-1\pmod{p^r}$ holds,
and this implies the claim.

For a given multiset $\fM_{p,r,a}(b)$, let $d$ be the least
positive integer congruent to $b$ modulo $p^r$, and let $\fM$ be
the multiset in $\Z$ defined by
$$
\fM=\big\langle\mathop{\underbrace{d,d,\ldots,d}}
\limits_{\text{$n$
copies}}~,\,\mathop{\underbrace{-d,-d,\ldots,-d}}
\limits_{\text{$p^r-2-n$ copies}}\big\rangle.
$$
Then $\fM$ and $\fM_{p,r,a}(b)$ are congruent modulo $p^r$.

Suppose first that $p=2$.  Since $\fM_{p,r,a}(b)=\fM_{p,r,a}(-b)$,
then replacing $b$ by $-b$ if necessary, we can assume that $d\le
2^r-d$. Hence, $d$ is a positive odd integer with $d\le 2^{r-1}$.
Also,
$$
n\equiv 2^{r-1}b^{-1}-1\equiv 2^{r-1}-1\pmod{2^r},
$$
where the second congruence follows from the fact that $b$ is odd;
in view of the minimality condition on $n$, it follows that
$n=2^{r-1}-1$.  Therefore, $\fM=\fN_{p,r}(\eta,c)$ with $\eta=n$
and $c=d$.

Now suppose that $p>2$. Since $\fM_{p,r,a}(b)=\fM_{p,r,a}(-b)$,
then replacing $b$ by $-b$ if necessary, we can assume that $n\le
p^r-2-n$. Let $\lambda$ be the least positive integer such that
$\lambda\equiv b^{-1}\pmod{p}$; then,
$$
n\equiv p^{r-1}b^{-1}-1\equiv p^{r-1}\lambda-1\pmod{p^r}.
$$
In view of the minimality condition on $n$ and the fact that
$n\le(p^r-2)/2$, it follows that
$n\in\{p^{r-1}\lambda-1:1\le\lambda\le(p-1)/2\}$. Also, defining
$\overline\lambda$ as in the statement of the lemma, we have
$$
d\equiv b\equiv\overline\lambda\pmod{p}.
$$
Since $1\le d\le p^r-1$, it follows that
$d\in\{\overline\lambda+p\mu:0\le\mu\le p^{r-1}-1\}$. Therefore,
$\fM=\fN_{p,r}(\eta,c)$ with $\eta=n$ and $c=d$.

To prove the uniqueness assertion, we must show that the multisets
$\fN_{p,r}(\eta,c)$ defined in the statement of the lemma are all
distinct modulo~$p^r$.  If $p^r=2$, then $\eta=0$, $c=1$, and
$\fN_{p,r}(\eta,c)=\varnothing$, so there is nothing to prove;
hence, we can assume that $p^r>2$. Now suppose that
\begin{equation}
\label{eq:congNN}
\fN_{p,r}(\eta_1,c_1)\equiv\fN_{p,r}(\eta_2,c_2)\pmod{p^r}.
\end{equation}
If $p=2$, we have $\eta_1=\eta_2=2^{r-1}-1$.  Also, since $c_j<
2^r-c_j$ for $j=1,2$ (note that the inequalities are strict since
$2^r\ge 4$), the congruence~\eqref{eq:congNN} implies that
$c_1\equiv c_2\pmod{2^r}$; as
$c_1,c_2\in\{1,3,\ldots,2^{r-1}-1\}$, this is possible only if
$c_1=c_2$.  If $p>2$, then the inequalities $\eta_j<p^r-2-\eta_j$,
$j=1,2$, and the congruence~\eqref{eq:congNN} together imply that
$\eta_1=\eta_2$ and $c_1\equiv c_2\pmod{p^r}$.  Since $1\le c_j\le
p^r-1$ for $j=1,2$, it follows that $c_1=c_2$.  This completes the
proof.
\end{proof}

\subsection{Algebraic results}

Let $G$ be a fixed \emph{nontrivial} cyclic group, and let $a$ be
an element of $G$ other than the identity.  Among the subgroups
$H< G$ that do not contain $a$, let $H(a)$ denote that subgroup
$H$ which has the greatest cardinality; note that $H(a)$ is
well-defined since every subgroup of a finite cyclic group is
determined uniquely by its cardinality.  Let
\begin{equation}
\label{eq:|G|} \big|G\big|=\prod_{j=1}^k
p_j^{\alpha_j}\qquad\text{and} \qquad {\rm
ord\/}_G(a)=\prod_{j=1}^s p_j^{\beta_j}
\end{equation}
be the prime factorizations of $\big|G\big|$ and ${\rm
ord\/}_G(a)$ (the order of $a$ in $G$).  Here, $p_1,\ldots,p_k$
are distinct primes, $s\le k$, and the integers $\alpha_j$ and
$\beta_j$ are positive. Using these data, we define:
$$
\r(G,a)=\min_{1\le j\le s}\big\{p_j^{\alpha_j-\beta_j+1}\big\}.
$$
Note that the constant $\paq$ defined by~\eqref{eq:paq_defn} is
equal to $\r(G,a)$ in the case that $G$ is the cyclic group
$(\Z/q\Z)^*$.

\begin{lemma}
\label{lem:G/H(a)} Let $G$ be a nontrivial cyclic group, and let
$a$ be an element of $G$ other than the identity. Then
$\r(G,a)\big|H(a)\big|=\big|G\big|$.
\end{lemma}

Thus, the definition of $H(a)$ given here is consistent with the
definition given in the introduction.

\begin{proof}
We begin by factoring $\big|G\big|$ and ${\rm ord\/}_G(a)$ as
in~\eqref{eq:|G|} above. By the \emph{Chinese Remainder Theorem},
we have
$$
G\cong\Z/\big|G\big|\Z\cong\prod_{j=1}^k \Z/p_j^{\alpha_j}\Z.
$$
Under this isomorphism, the element $a\in G$ can be identified
with an ordered $k$-tuple:
$$
\big(p_1^{\alpha_1-\beta_1}a_1,\ldots,
p_k^{\alpha_k-\beta_k}a_k\big)\in\prod_{j=1}^k
\Z/p_j^{\alpha_j}\Z,
$$
where each $a_j$ is an integer in the range
     $1\le a_j\le p_j^{\beta_j}$, and $p_j\nmid a_j$. Replacing
$a$ by one of its automorphic images $\widetilde a\in G$, if
necessary, we can assume without loss of generality that $a_j=1$
for $j=1,\ldots,k$. Indeed, if a subset $S\subset G$ is {\it
characteristic\/} (that is, invariant under all automorphisms of
$G$), then $S$ does not contain $a$ if and only if $S$ does
contain any automorphic image of $a$.  Since $H(a)$ is
characteristic, it follows that $H(a)=H(\widetilde a)$ for every
automorphic image $\widetilde a$ of $a$.

Now let $K$ be a subgroup of $G$ that does not contain $a$, and
suppose that $\big|K\big|=\prod_{j=1}^kp_j^{\gamma_j}$ for some
nonnegative integers $\gamma_j$.  Then,
$$
K\cong\prod_{j=1}^k p_j^{\alpha_j-\gamma_j}\Z/p_j^{\alpha_j}\Z
\subseteq\prod_{j=1}^k \Z/p_j^{\alpha_j}\Z.
$$
The condition that $a\not\in K$ is equivalent to the existence of
an index $j$ such that $\alpha_j-\gamma_j>\alpha_j-\beta_j$; that
is, $\gamma_j<\beta_j$. In particular, $\beta_j>0$, and therefore
$1\le j\le s$.  If $K$ is \emph{maximal} among the subgroups of
$G$ which do not contain $a$, it must be the case that
$\gamma_j=\beta_j-1$ and $\gamma_i=\alpha_i$ for all $i\ne j$;
consequently, $\big|K\big|=\big|G\big|/p_j^{\alpha_j-\beta_j+1}$.
Finally, since $H(a)$ has the largest cardinality of all such
subgroups $K$, it is clear that
$\big|H(a)\big|=\big|G\big|/p_j^{\alpha_j-\beta_j+1}$, where $j$
is the only index for which
$$
p_j^{\alpha_j-\beta_j+1}=\min_{1\le i\le
s}\big\{p_i^{\alpha_i-\beta_i+1}\big\}=\r(G,a),
$$
and this completes the proof.
\end{proof}

\begin{lemma}
\label{lem:G/H(a)two} Let $a$ be a fixed element of $G=(\Z/q\Z)^*$
other than the identity, and suppose that $\paq$ is the prime
power $p^r$. Put $G_r=\Z/p^r\Z$. Then there exists a generator $g$
of the cyclic group $G$ such that the map $n\mapsto g^n$ defines a
group isomorphism $\phi_g:G_r\to G/H(a)$ which maps the congruence
class $p^{r-1}\pmod{p^r}$ to the coset $aH(a)$.
\end{lemma}

\begin{proof}
First, let $g$ be an arbitrary generator of $G$. Since every
subgroup of $G$ is determined uniquely by its cardinality, it
follows from Lemma~\ref{lem:G/H(a)} that $H(a)$ is the subgroup of
$G$ generated by $g^{\paq}=g^{p^r}$.  Then, it is easy to see that
the map $n\mapsto g^n$ defines a group isomorphism $\phi_g:G_r\to
G/H(a)$.

Let $\psi_g:G\to G_r$ be the homomorphism defined via the
composition:
$$
G\to G/H(a)\mathop{\,\longrightarrow\,}\limits^{\phi_g^{-1}} G_r.
$$
Since $a\not\in H(a)$, the element $\overline a=\psi_g(a)$ is not
the identity in $G_r$. On the other hand, $\overline a$ is
contained in \emph{every} subgroup $K$ of $G_r$, for otherwise the
preimage $\psi_g^{-1}(K)$ would be a subgroup of $G$ which
properly contains $H(a)$ and such that $a\not\in\psi_g^{-1}(K)$,
contradicting the maximality of $H(a)$. In particular, $\overline
a$ lies in the subgroup $K$ generated in $G_r$ by the congruence
class $p^{r-1}\pmod{p^r}$. Thus, $\psi_g(a)=b p^{r-1}\pmod{p^r}$
for some integer $b$ with $p\nmid b$. Replacing $g$ by the
generator $g^b$, the result follows immediately.
\end{proof}

The following technical lemma, used in the proof of
Theorem~\ref{thm:hard} below, combines the preceding two lemmas
with the combinatorial results of the previous section.

\begin{lemma}
\label{lem:technical} Let $a$ be a fixed element of $G=(\Z/q\Z)^*$
other than the identity. Write $\paq=p^r$, and put $G_r=\Z/p^r\Z$.
Let $g$ be a generator of $G$ with the property described in
Lemma~\ref{lem:G/H(a)two}, and let $\psi_g:G\to G_r$ be the
homomorphism defined in the proof of that lemma.

Suppose that $\fM=\langle a_1,\ldots,a_k\rangle$ is a multiset in
$G$ with the property:
\begin{equation}
\label{eq:zeroprod1} \prod_{j\in S} a_j\ne a \quad\text{for every
subset~}S\subseteq \{1,2,\ldots,k\}.
\end{equation}
Let $\fH$ be the multiset consisting of the elements $a_j\in\fM$
that occur with multiplicity at least $q-2$, and let $H$ be the
subgroup of $G$ generated by the elements of $\fH$.  Finally, let
$\fK$ be the multiset consisting of those elements of $\fM$ which
do not lie in $H$. Then:
\begin{itemize}
\item[$(i)$] $|\fK|\le(q-1)(q-3)$;

\item[$(ii)$] $|H|\le|H(a)|$, and equality holds if and only if
$H=H(a)$;

\item[$(iii)$] If $H=H(a)$, then $|\fK|\le\paq-2$;

\item[$(iv)$] Suppose that $H=H(a)$ and $|\fK|=\paq-2$. Then
$\psi_g(\fK)$ is a multiset $\fN_{p,r}(\eta,c)$ of the type
considered in Lemma~\ref{lem:listmultisets}.
\end{itemize}
\end{lemma}

\begin{proof}
The assertion $(i)$ is trivial, since $|G|=q-1$ and every element
of $\fK$ occurs with multiplicity at most $q-3$.

Let $b_1,\ldots,b_s$ be the distinct elements that occur in the
multiset $\fH$, and let $m_1,\ldots,m_s$ be the respective
multiplicities. Since every element of $H$ can be expressed as a
product $\prod_{i=1}^sb_i^{\nu_i}$, where $0\le\nu_i\le q-2\le
m_i$ for $i=1,\ldots,s$, it follows that every element of $H$ is a
product of the form $\prod_{j\in S}a_j$ for some subset
$S\subseteq\{1,\ldots,k\}$.  Using~\eqref{eq:zeroprod1}, we see
that $a\not\in H$, hence $(ii)$ follows immediately from the
definition of $H(a)$ and the fact that every subgroup of $G$ is
determined uniquely by its cardinality.

{From} now on, we assume $H=H(a)$. Write $\fK=\langle
k_1,\ldots,k_t\rangle$, and observe that
\begin{equation}
\label{eq:zeroprod2} \prod_{i\in T} k_i\not\in aH(a)
\quad\text{for every subset~}T\subseteq \{1,2,\ldots,t\}.
\end{equation}
Indeed, assuming that $\prod_{i\in T}k_i=ah^{-1}$ for some $h\in
H(a)$, the argument above shows that $h=\prod_{j\in S}a_j$ for
some subset $S\subseteq\{1,\ldots,k\}$, and as $\fK\subseteq\fM$,
it follows that $\prod_{i\in T}k_i=\prod_{j\in R}a_j$ for another
subset $R\subseteq\{1,\ldots,k\}$.  Clearly, $R\cap S=\varnothing$
since $\fK\cap\fH=\varnothing$; therefore, $\prod_{j\in R\cup
S}a_j=a$, which contradicts~\eqref{eq:zeroprod1}.

Let $\overline\fK=\langle \overline k_1,\ldots,\overline
k_t\rangle$ be the image of $\fK$ under the map $\psi_g$, that is
$\overline\fK=\psi_g(\fK)$, and put $\overline a=\psi_g(a)$.
Using~\eqref{eq:zeroprod2}, we deduce that
$$
\prod_{i\in T}\, \overline k_i\ne\overline a \quad\text{for every
subset~}T\subseteq \{1,2,\ldots,t\}.
$$
Therefore, Lemma~\ref{lem:comb2} immediately implies that
$$
\big|\fK\,\big|=\big|\,\overline\fK\,\big|=t\le\kappa(G_r,\overline
a)= p^r-2=\paq-2,
$$
which proves $(iii)$. In the case that $|\fK\,|=\paq-2$, we can
apply Lemmas~\ref{lem:extremecase} and~\ref{lem:listmultisets} to
conclude that $\overline\fK=\fN_{p,r}(\eta,c)$ for a unique choice
of $\eta$ and $c$, which proves $(iv)$.
\end{proof}

\subsection{Analytic results}

For the proofs of Theorems~\ref{thm:easy} and~\ref{thm:hard}, we
need a variant of the classical result of Landau~\cite{La}:
\begin{equation}
\label{eq:LandauResult} \big|\{n\le x~:~\Omega(n)=k\}\big|
\sim\big|\{n\le x~:~\omega(n)=k\}\big|
\sim\frac{x(\log_2x)^{k-1}}{(k-1)!\log x},
\end{equation}
where $k\ge 1$ is a fixed integer, and $\Omega(n)$ and $\omega(n)$
denote the total number of prime factors of $n$ counted with and
without multiplicity, respectively. Specifically, we need an
estimate for number of positive integers $n\le x$ with
$\Omega(n)=k$ and such that every prime factor of $n$ lies in a
prescribed subset of the residue classes modulo $m$.

In this section the implied constants, frequently without explicit
mention, may depend not only on $m$ but on $k$ and on various
other parameters; virtually everything but $x$ is fair game.

For given $m$ let $\cA$ be a nonempty subset of $(\Z/m\Z)^*$.
Define
$$
\cQ(m,\cA)=\{n\ge 1~:~p\mid n\,\Rightarrow\,p\equiv a\pmod
m\text{~for some~}a\in\cA\}.
$$
For each $k$ define $\cQ_k(m,\cA)$ to be the set of positive
integers $n$ in $\cQ(m,\cA)$ for which $\Omega(n)=k$.

\begin{lemma}
\label{lem:anal1} Let $m$ and $k$ be fixed positive integers and
$\cA$ a nonempty subset of $(\Z/m\Z)^*$. For real $x\ge 1$, let
$\cQ_k(x;m,\cA)$ be the number of positive integers $n\le x$ in
the set $\cQ_k(m,\cA)$. Then,
$$
\cQ_k(x;m,\cA)=(1+o(1))\(\frac{|\cA|}{\varphi(m)}\)^k
\frac{x(\log_2 x)^{k-1}}{(k-1)!\,\log x}\ .
$$
\end{lemma}

\begin{proof} For the proof we may follow an argument given in
Section~9.4 in the book by Nathanson~\cite{Nat0}.

Let $\cP$ be the set of primes $p$ such that $p\equiv a\pmod m$
for some $a\in\cA$, let $\cP^k$ be the set of ordered $k$-tuples
of primes in $\cP$, and for every positive integer $n$, let
$$
r_k(n)=\big|\big\{(p_1,\ldots,p_k)\in\cP^k~:~p_1\cdots
p_k=n\big\}\big|.
$$
For any real number $x\ge 1$, put
\begin{eqnarray*}
f_k(x)&=&\sum_{n\le x}r_k(n)
=\sum_{\substack{p_1\cdots p_k\le x\\(p_1,\ldots,p_k)\in\cP^k}}1,\\
g_k(x)&=&\sum_{n\le x}\frac{r_k(n)}{n} =\sum_{\substack{p_1\cdots
p_k\le x\\(p_1,\ldots,p_k)\in\cP^k}}
\frac{1}{p_1\cdots p_k},\\
h_k(x)&=&\sum_{n\le x}r_k(n)\log n =\sum_{\substack{p_1\cdots
p_k\le x\\(p_1,\ldots,p_k)\in\cP^k}} \log(p_1\cdots p_k).
\end{eqnarray*}
Note that, for every $k\ge 1$, the relations
\begin{equation}
\label{eq:relt1} g_{k+1}(x)=\sum_{p\le
x,~p\in\cA}\frac{g_k(x/p)}{p},
\end{equation}
and
\begin{equation}
\label{eq:relt2} k\,h_{k+1}(x)=(k+1)\sum_{p\le x,~p\in\cA}h_k(x/p)
\end{equation}
follow easily from the above definitions. Finally, let
$\cQ_k^{\flat}(m,\cA)$ denote the set of all {\it squarefree\/}
elements of $\cQ_k(m,\cA)$ and let $\cQ_k^{\flat}(x;m,\cA)$ count
the number of these up to $x$. For these, we of course have
$\Omega(n)=\omega(n)=k$.

The following properties of $r_k(n)$ are immediate:
\begin{itemize}
\item $0\le r_k(n)\le k!$ for all $n\ge 1$;

\item $r_k(n)>0~\Longleftrightarrow~n\in\cQ_k(m,\cA)$;

\item $r_k(n)=k!~\Longleftrightarrow~n\in\cQ_k^{\flat}(m,\cA)$.
\end{itemize}
Consequently,
\begin{equation}
\label{eq:Land1} f_k(x)=\sum_{n\le x}r_k(n)\le
k!\sum_{\substack{n\le x\\r_k(n)>0}}1 =k!\,\cQ_k(x;m,\cA),
\end{equation}
and
\begin{equation}
\label{eq:Land2} f_k(x)=\sum_{n\le x}r_k(n)\ge
k!\sum_{\substack{n\le x\\r_k(n)=k!}}1
=k!\,\cQ_k^{\flat}(x;m,\cA).
\end{equation}
If $n\in\cQ_k(m,\cA)\setminus\cQ_k^{\flat}(m,\cA)$, then
$\omega(n)<\Omega(n)=k$, and therefore,
$$
\cQ_k(x;m,\cA) - \cQ_k^{\flat}(x;m,\cA)\le\big|\{n\le
x~:~\omega(n)<k\}\big| \ll\frac{x(\log_2x)^{k-2}}{\log x},
$$
where we have used~\eqref{eq:LandauResult} in the last step.
Hence, from~\eqref{eq:Land1} and~\eqref{eq:Land2} we deduce that
$$
\cQ_k(x;m,\cA)=\frac{f_k(x)}{k!}+O\(\frac{x(\log_2x)^{k-2}}{\log
x}\).
$$
To prove the theorem, it therefore suffices to establish the
estimate:
\begin{equation}
\label{eq:desired_est1} f_k(x)=(1+o(1))\,\frac{C^kkx(\log_2
x)^{k-1}}{\log x},
\end{equation}
where, for brevity, we have put $C=|\cA|/\varphi(m)$. As it is
clear that $f_k(x)=O(x)$, by partial summation we have
$$
h_k(x)=\sum_{n\le x}r_k(n)\log n=f_k(x)\log x-\int_1^x
\frac{f_k(t)}{t}\,dt =f_k(x)\log x+O(x),
$$
and thus~\eqref{eq:desired_est1} follows immediately from the
estimate:
\begin{equation}
\label{eq:desired_est2} h_k(x)=(1+o(1))\,C^kkx(\log_2x)^{k-1}.
\end{equation}

Using the prime number theorem for arithmetic progressions we have
\begin{equation}
\label{eq:h1x} h_1(x)=\sum_{p\le x,~p\in\cP}\log p =(1+o(1))\,C x.
\end{equation}
In particular, this yields~\eqref{eq:desired_est2} in the special
case $k=1$. By the analogue for arithmetic progressions of the
theorem of Mertens, or by partial summation from the previous
formula, we also have
$$
g_1(x)=\sum_{p\le x,~p\in\cP}\frac{1}{p} =(1+o(1))\,C\log_2x.
$$
The latter estimate implies that
$$
g_1(x^{1/k})=(1+o(1))\,C\log_2(x^{1/k})=(1+o(1))\,C\log_2x.
$$
Thus, from the trivial inequalities
$$
g_1(x^{1/k})^k\le g_k(x)\le g_1(x)^k,
$$
we see that
\begin{equation}
\label{eq:gkest} g_k(x)=(1+o(1))\,C^k(\log_2x)^k.
\end{equation}

Now, for $k\ge 1$, define
$$
F_k(x)=h_k(x)-C kxg_{k-1}(x),
$$
where we have put $g_0(x)=1$ for all $x\ge 1$. We claim that the
bound
\begin{equation}
\label{eq:Fkissmall} F_k(x)=o(x(\log_2x)^{k-1})
\end{equation}
holds for each fixed $k\ge 1$. Observe that this statement implies
the desired result; indeed, if~\eqref{eq:Fkissmall} holds for some
integer $k\ge 2$, then in view of~\eqref{eq:gkest}, we have
$$
h_k(x)=C kxg_{k-1}(x)+F_k(x) =C^k kx(\log_2x)^{k-1}
+o(x(\log_2x)^{k-1}),
$$
which gives~\eqref{eq:desired_est2}.

To prove~\eqref{eq:Fkissmall}, we use induction on the parameter
$k$. The initial case $k=1$ follows immediately
from~\eqref{eq:h1x} and the fact that $g_0(x)=1$.  Now suppose
that~\eqref{eq:Fkissmall} holds for some integer $k\ge 1$. By
relations~\eqref{eq:relt1} and~\eqref{eq:relt2}, we have
\begin{eqnarray}
k\,F_{k+1}(x)&=&k\,h_{k+1}(x)-C k(k+1)xg_k(x)\nonumber\\
&=&(k+1)\sum_{p\le x,~p\in\cA}h_k(x/p)
+C k(k+1)x\sum_{p\le x,~p\in\cA}\frac{g_{k-1}(x/p)}{p}\nonumber\\
&=&(k+1)\sum_{p\le x,~p\in\cA}\Big(h_k(x/p)
+C k(x/p)g_{k-1}(x/p)\Big)\nonumber\\
\label{eq:blue} &=&(k+1)\sum_{p\le x,~p\in\cA}F_k(x/p).
\end{eqnarray}
For fixed $\eps>0$, there is a constant $x_0=x_0(\eps)$ such that
$$
F_k(x/p)\le\frac{\eps x(\log_2(x/p))^{k-1}}{p} \le\frac{\eps
x(\log_2x)^{k-1}}{p},
$$
whenever $x/p\ge x_0$.  On the other hand, $F_k(x/p)=O_\eps(1)$ if
$x/p<x_0$. Consequently,
\begin{eqnarray*}
\sum_{p\le x,~p\in\cA}F_k(x/p)&=&\sum_{p\le
x/x_0,~p\in\cA}F_k(x/p)+
\sum_{x/x_0<p\le x,~p\in\cA}F_k(x/p)\\
&\le&\eps x(\log_2x)^{k-1}\sum_{p\le x/x_0}\frac{1}{p}+
O_\eps\(\sum_{x/x_0<p\le x}1\)\\
&=&(\eps+o(1))x(\log_2x)^k+O_\eps(x/\log x).
\end{eqnarray*}
Combining this estimate with~\eqref{eq:blue}, it follows that for
every $\eps>0$, there is a constant $x_1=x_1(\eps)$ such that the
inequality
$$
F_{k+1}(x)\le 2\eps(1+1/k)x(\log_2x)^k
$$
holds whenever $x\ge x_1$; it follows that~\eqref{eq:Fkissmall}
holds with $k$ replaced by $k+1$.  This completes the induction
and finishes the proof of the lemma.
\end{proof}

We need to count, in addition to the integers in $\cQ_k(m,\cA)$,
the same integers without the restriction on $\Omega(n)$, that is
those in the set $\cQ(m,\cA)$. For this we shall use the following
result of Wirsing~\cite{Wirs}:

\begin{lemma}
\label{lem:wirs} Suppose that $f$ is a fixed real-valued
multiplicative function with the following properties:
\begin{itemize}
\item[$(i)$] For every natural number $n$, $f(n)\ge 0$;

\item[$(ii)$] For some constants $c_1,c_2$ with $c_2<2$, the bound
$f(p^\nu)\le c_1\,c_2^\nu$ holds for all primes $p$ and integers
$\nu\ge 2$;

\item[$(iii)$] There exists a constant $\tau > 0$ such that
$$
\sum_{p \le x} f(p)  = \(\tau + o(1) \) \frac{x}{\log x}.
$$
\end{itemize}
Then,
$$
\sum_{n \le x} f(n) = \(\frac{1} {e^{\gamma\tau}\Gamma(\tau)}
+o(1)\) \frac{x}{\log x}\prod_{p\leq x}\,  \sum_{\nu=0}^\infty
\frac{f(p^\nu)}{p^\nu},
$$
where $\gamma$ is the Euler-Mascheroni constant and $\Gamma(s)$ is
the gamma function.
\end{lemma}

The classical result of Mertens that
$$
\prod_{p\le x}\(1-\frac{1}{p}\)
     = e^{-\gamma}(\log x)^{-1}+O\((\log x)^{-2}\),
$$
has been generalized in the paper of Williams~\cite{Williams} (see
also~\cite{Uchiyama}), which gives a similar estimate when the
product above is restricted to primes lying in a fixed arithmetic
progression. To state this result we first recall some notation
from~\cite{Williams}. Let $m$ be a positive integer and let $\chi$
be a non-principal Dirichlet character modulo~$m$. Let $L(s,\chi)$
be the corresponding $L$--function and define the Dirichlet series
$$
K(s,\chi)=\sum_{n=1}^\infty\frac{k_\chi(n)}{n^s}
     =\prod_p\(1-\frac{k_\chi(p)}{p^s}\)^{-1},
$$
where $k_\chi(n)$ is the completely multiplicative function whose
value at the prime~$p$ is given by
$$
k_\chi(p)=p\(1-\(1-\frac{\chi(p)}{p}\)
     \(1-\frac{1}{p}\)^{-\chi(p)}\).
$$
The main result of~\cite{Williams} is the following:

\begin{lemma}
\label{lem:Lem} Let $a$ and $m\ge 1$ be coprime integers. Then,
\begin{equation}
\label{eq:needO} \prod_{\substack{p\le x\\p\equiv a\pmod
m}}\(1-\frac{1}{p}\)=\frac{\varpi(a,m)}{(\log
x)^{1/\varphi(m)}}+O\(\frac{1}{(\log x)^{1+1/\varphi(m)}}\),
\end{equation}
where
\begin{equation}
\label{eq:define-varpi}
\varpi(a,m)=\(e^{-\gamma}\frac{m}{\varphi(m)}
     \prod_{\chi\ne\chi_0}\(\frac{K(1,\chi)}{L(1,\chi)}
     \)^{\overline\chi(a)}\)^{1/\varphi(m)}.
\end{equation}
\end{lemma}

We are now ready to count the integers in $\cQ(m,\cA)$. Recall
that these are just the integers all of whose prime factors lie in
the set $\cA$. For real $x\ge 1$ let $\cQ(x;m,\cA)$ denote the
number of such integers $n\le x$.
\begin{lemma}
\label{lem:A}  Let $m$ be a fixed positive integer and $\cA$ a
nonempty subset of $(\Z/m\Z)^*$. Then,
$$
\cQ(x;m,\cA)=(1+o(1))\,\vartheta(m,\cA) \frac{x}{(\log
x)^{1-|A|/\varphi(m)}},
$$
where
\begin{equation}
\label{eq:define-theta} \vartheta(m,\cA)= \frac{e^{-\gamma
|A|/\varphi(m)}}{\Gamma\big(|A|/\varphi(m)\big)} \prod_{a\in
\cA}\varpi(a,m)^{-1},
\end{equation}
with the constants $\varpi(a,m)$ defined as in
Lemma~\ref{lem:Lem}.
\end{lemma}

\begin{proof}
This follows immediately by applying Lemma~\ref{lem:wirs} to the
multiplicative function $f$ which is defined on prime powers by
$$
f(p^\nu)=\left\{
\begin{array}{ll}
        1 &\quad\hbox{if $p\equiv a\pmod m$ for some $a\in\cA$;} \\
        0 &\quad\hbox{otherwise;} \\
\end{array}
\right.
$$
making use of the estimates of Lemma~\ref{lem:anal1} (with $k=1$)
and of Lemma~\ref{lem:Lem}.
\end{proof}

The next lemma evaluates $\vartheta(m,\cA)$ in the special case
$m=q$, $\cA=H(a)$.

\begin{lemma}
\label{lem:theta2V} We have
$$
\vartheta(q,H(a))
=\frac{e^{-\gamma/\paq}(1-q^{-1})^{1/\paq}}{\Gamma(1/\paq)}\,\vaq,
$$
where $\paq$ and $\vaq$ are given by~\eqref{eq:paq_defn}
and~\eqref{eq:vaq_defn}, respectively.
\end{lemma}

\begin{proof}
By the definitions~\eqref{eq:define-varpi}
and~\eqref{eq:define-theta}, we have
$$
\vartheta(q,H(a))=\frac{e^{-\gamma/\paq}}{\Gamma\big(1/\paq\big)}
\prod_{h\in H(a)}\varpi(h,q)^{-1},
$$
where
$$
\varpi(h,q)=\(e^{-\gamma}\frac{q}{q-1}
     \prod_{\chi\ne\chi_0}\(\frac{K(1,\chi)}{L(1,\chi)}
     \)^{\overline\chi(h)}\)^{1/(q-1)}.
$$
{From} the orthogonality relation
$$
\sum_{h\in H(a)}\overline\chi(h)=\left\{
\begin{array}{ll}
      |H(a)| & \quad\hbox{if $\chi\big|_{H(a)}=1$,} \\
      \quad 0 & \quad\hbox{otherwise,} \\
\end{array}
\right.
$$
it follows that
$$
\vartheta(q,H(a))=\frac{(1-q^{-1})^{1/\paq}}{\Gamma(1/\paq)}
\prod_{\substack{\chi\ne\chi_0\\\chi|_{H(a)}=1}}
\(\frac{K(1,\chi)}{L(1,\chi)}\)^{-1/\paq}.
$$
By $(3.2)$ of~\cite{Williams}, we have
$$
\frac{K(1,\chi)}{L(1,\chi)}=\lim_{y\to\infty}\prod_{p\le
y}\(1-\frac{1}{p}\)^{\chi(p)}.
$$
Therefore, in view of the relation
$$
\sum_{\substack{\chi\ne\chi_0\\\chi|_{H(a)}=1}}
\chi(p)=\left\{%
\begin{array}{ll}
      \paq-1 & \quad\hbox{if $p\!\!\!\pmod q\in H(a)$,} \\
      \quad -1 & \quad\hbox{otherwise,} \\
\end{array}%
\right.
$$
and the Mertens' formula
$$
\prod_{p\le y}\(1-\frac{1}{p}\)
   =(1+o(1))\frac{e^{-\gamma}}{\log y},
$$
we derive that
\begin{equation*}
\begin{split}
\vartheta(q,&H(a))=\frac{(1-q^{-1})^{1/\paq}}{\Gamma(1/\paq)}
\Biggl(\lim_{y\to\infty}\prod_{p\le
y}\prod_{\substack{\chi\ne\chi_0\\\chi|_{H(a)}=1}}
\(1-\frac{1}{p}\)^{\chi(p)}\Biggr)^{-1/\paq}\\
&=\frac{(1-q^{-1})^{1/\paq}}{\Gamma(1/\paq)}
\lim_{y\to\infty}\Biggl(\prod_{p\le
y}\(1-\frac{1}{p}\)^{-1}\!\!\!\!\!\!\prod_{\substack{p\le
y\\p\!\!\!\pmod q\in H(a)}}
\!\!\!\!\!\!\(1-\frac{1}{p}\)^{\paq}\Biggr)^{-1/\paq}\\
&=\frac{e^{-\gamma/\paq}(1-q^{-1})^{1/\paq}}{\Gamma(1/\paq)}
\lim_{y\to\infty}\Biggl((\log y)^{-1/\paq}\!\!\!\!\!\!
   \prod_{\substack{p\le y\\p\!\!\!\pmod q\in
H(a)}}\!\!\!\!\!\!\(1-\frac{1}{p}\)^{-1}\Biggr).
\end{split}
\end{equation*}
Inserting the definition~\eqref{eq:vaq_defn}, we finish the proof.
\end{proof}

\begin{lemma}
\label{lem:sumSxxx} Let $a$ be a nonnegative integer, $b$ a real
number in the half-open interval $(0,1]$, $c$ a nonnegative real
number, and $K$ a positive real number. Let $\cS$ be a set of
positive integers, and for $x\ge 1$ put
$$
\cS(x)=|\{n\le x~:~n\in\cS\}|\ .
$$
Finally, suppose that the following estimate holds as
$x\to\infty$:
$$
S(x)=(1+o(1))\frac{Kx(\log_2x)^a}{(\log x)^b}.
$$
Then,
$$
\sum_{\substack{h\in\cS\\h\le x^{1/2}}}\frac{1}{h\,(\log(x/h))^c}=\left\{%
\begin{array}{ll}
      \displaystyle
      O\(\frac{(\log_2x)^a}{(\log x)^{b+c-1}}\)
       &\quad\hbox{if $b\in(0,1)$;} \\ \\
      \displaystyle
      (1+o(1))\frac{K}{a+1}\frac{(\log_2x)^{a+1}}{(\log x)^c}
       &\quad\hbox{if $b=1$.} \\
\end{array}%
\right.
$$
\end{lemma}

\begin{proof}
Since
$$
\sum_{h\le\log_2x}\frac{1}{h\,(\log(x/h))^c}
=(1+o(1))\frac{\log_3x}{(\log x)^c},
$$
we have
$$
\sum_{\substack{h\in\cS\\h\le x^{1/2}}}\frac{1}{h\,(\log(x/h))^c}=
\sum_{\substack{h\in\cS\\\log_2x<h\le x^{1/2}}}
\frac{1}{h\,(\log(x/h))^c} +O\(\frac{\log_3x}{(\log x)^c}\).
$$
Since the estimate
$$
S(t)=(K+o(1))\frac{t(\log_2t)^a}{(\log t)^b}
$$
holds uniformly for all $t\ge\log_2x$, by partial summation we
deduce that
$$
\sum_{\substack{h\in\cS\\\log_2x<h\le x^{1/2}}}
\frac{1}{h\,(\log(x/h))^c}=\int_{\log_2x}^{x^{1/2}}
\frac{dS(t)}{t\,(\log(x/t))^c}=(K+o(1))(J_1+J_2-J_3),
$$
where
\begin{eqnarray*}
J_1&=&\biggl[\frac{(\log_2t)^a}{(\log t)^b(\log(x/t))^c}
\biggl]_{\log_2x}^{x^{1/2}},\\
J_2&=&\int_{\log_2x}^{x^{1/2}}\frac{(\log_2t)^a}{(\log
t)^b(\log(x/t))^c}\,\frac{dt}{t},\\
J_3&=&c\int_{\log_2x}^{x^{1/2}}\frac{(\log_2t)^a}{(\log
t)^b(\log(x/t))^{c+1}}\,\frac{dt}{t}.
\end{eqnarray*}
Clearly,
$$
J_1\ll\frac{(\log_4x)^a}{(\log_3x)^b(\log x)^c}\mand
J_3\ll\frac{J_2}{\log x}.
$$
Making the change of variables $t=x^s$ in the integral $J_2$, it
follows that
$$
J_2=\frac{1}{(\log x)^{b+c-1}}\int_{u_0}^{1/2}\frac{(\log_2x+\log
s)^a}{s^b\,(1-s)^c}\,ds,
$$
where $u_0=(\log_3x)/\log x$.  To complete the proof, it suffices
to show that
$$
\int_{u_0}^{1/2}\frac{(\log_2x+\log s)^a}{s^b\,(1-s)^c}\,ds
=\left\{
\begin{array}{ll}
      \displaystyle
      O\big((\log_2x)^a\big)
       &\quad\hbox{if $b\in(0,1)$;} \\ \\
      \displaystyle
      (1+o(1))\frac{(\log_2x)^{a+1}}{a+1}
       &\quad\hbox{if $b=1$.} \\
\end{array}
\right.
$$

First, we discuss the integral over $u_0\le s\le u_1$, where
$u_1=(\log_3x)^{-1}$. If $b\in(0,1)$, then
\begin{equation*}
\begin{split}
&\int_{u_0}^{u_1}\frac{(\log_2x+\log
s)^a}{s^b\,(1-s)^c}\,ds=(1+o(1))
\int_{u_0}^{u_1}\frac{(\log_2x+\log s)^a}{s^b}\,ds\\
&\qquad=(1+o(1))\biggl[\frac{s^{1-b}}{1-b}
\sum_{j=0}^a\frac{1}{(b-1)^j}\frac{a!}{(a-j)!}
   \,(\log_2x+\log s)^{a-j}\biggl]_{u_0}^{u_1}\\
&\qquad=(1+o(1))\frac{(\log_2x)^a}{(1-b)(\log_3x)^{1-b}},
\end{split}
\end{equation*}
and for $b=1$, we have
\begin{equation*}
\begin{split}
&\int_{u_0}^{u_1}\frac{(\log_2x+\log
s)^a}{s^b\,(1-s)^c}\,ds=(1+o(1))
\int_{u_0}^{u_1}\frac{(\log_2x+\log s)^a}{s}\,ds\\
&\qquad=(1+o(1))\biggl[\frac{(\log_2x+\log s)^{a+1}}{a+1}
\biggl]_{u_0}^{u_1}=(1+o(1))\frac{(\log_2x)^{a+1}}{a+1}.
\end{split}
\end{equation*}
Next, we consider the integral over $u_1\le s\le 1/2$. If
$b\in(0,1)$, then
\begin{equation*}
\begin{split}
\int_{u_1}^{1/2}\frac{(\log_2x+\log s)^a}{s^b\,(1-s)^c}\,ds
&=(1+o(1))(\log_2x)^a\int_{u_1}^{1/2}\frac{ds}{s^b\,(1-s)^c}\\
&\ll(\log_2x)^a,
\end{split}
\end{equation*}
and for $b=1$, we have
\begin{equation*}
\begin{split}
\int_{u_1}^{1/2}\frac{(\log_2x+\log s)^a}{s^b\,(1-s)^c}\,ds
&=(1+o(1))(\log_2x)^a\int_{u_1}^{1/2}\frac{ds}{s\,(1-s)^c}\\
&\ll(\log_2x)^a\log_3x.
\end{split}
\end{equation*}
Combining the preceding results, we obtain the stated estimates.
\end{proof}

\begin{lemma}
\label{lem:SSSSS} For $j=1,2$, let $\cS_j$ be a set of positive
integers and, for $x\ge 1,$ put
$$
\cS_j(x)=|\{n\le x~:~n\in\cS_j\}|\ .
$$
Suppose that, for $j=1, 2$,
$$
\cS_j(x)=(1+o(1))\frac{K_j\,x(\log_2x)^{a_j}}{(\log x)^{b_j}}
$$
where $a_1,a_2$ are nonnegative integers, $K_1,K_2>0$, and
$0<b_1<1$, $b_2=1$. Let $\cS(x)$ be the number of ordered pairs
$(h_1,h_2)\in\cS_1\times\cS_2$ such that $h_1h_2\le x$. Then the
following estimate holds:
$$
\cS(x)=(1+o(1))\,\frac{K_1K_2}{a_2+1}
\,\frac{x(\log_2x)^{a_1+a_2+1}}{(\log x)^{b_1}}.
$$
\end{lemma}

\begin{proof}
Observe that
$$
\cS(x)=\sum_{\substack{h_1\in\cS_1\\h_1\le x^{1/2}}}\cS_2(x/h_1)
   +\sum_{\substack{h_2\in\cS_2\\h_2\le x^{1/2}}}\cS_1(x/h_2)
   -\cS_1(x^{1/2})\,\cS_2(x^{1/2}).
$$
Uniformly for $h_1\le x^{1/2}$, we have
$$
\cS_2(x/h_1)=(1+o(1))
\,\frac{K_2\,x(\log_2(x/h_1))^{a_2}}{h_1\log(x/h_1)};
$$
thus Lemma~\ref{lem:sumSxxx} implies that
$$
\sum_{\substack{h_1\in\cS_1\\h_1\le x^{1/2}}}\cS_2(x/h_1)\ll
x(\log_2x)^{a_2}\sum_{h_1\in\cS_1(x^{1/2})}
\frac{1}{h_1\log(x/h_1)}\ll\frac{x(\log_2x)^{a_1+a_2}}{(\log
x)^{b_1}}.
$$
Similarly,
\begin{equation*}
\begin{split}
\sum_{\substack{h_2\in\cS_2\\h_2\le x^{1/2}}}\cS_1(x/h_2)
&=(1+o(1))\,K_1\,x(\log_2x)^{a_1}
\sum_{h_2\in\cS_2(x^{1/2})}\frac{1}{h_2(\log(x/h_2))^{b_1}}\\
&=(1+o(1))\,\frac{K_1K_2}{a_2+1}
\,\frac{x(\log_2x)^{a_1+a_2+1}}{(\log x)^{b_1}},
\end{split}
\end{equation*}
where we have again used Lemma~\ref{lem:sumSxxx}. Since
$$
\cS_1(x^{1/2})\,\cS_2(x^{1/2})\ll\frac{x(\log_2x)^{a_1+a_2}}{(\log
x)^{b_1+1}},
$$
the result follows.
\end{proof}

\section{Proofs of the theorems}

In this section we frequently use the notation $S(x)$ for the
number of positive integers $n\le x$ in the set $S$.

\subsection{Proof of Theorem~\ref{thm:easy}}
\label{sec:thm1}

Fix the prime $q\ge 3$, and write $\cN$ and $\cN(x)$ respectively
for $\cN(q,1)$, and $\cN(x;q,1)$. Let $\cN^*$ be the set of
integers $n\in\cN$ that are not divisible by $q$. Then $\cN^*$ can
be expressed as a disjoint union $\cN_1\cup\cN_2$, where $\cN_1$
is the set of integers $n\in\cN^*$ with $\Omega(n)\le q-3$, and
$\cN_2=\cN^\setminus\cN_1$.

Since $\cN_1$ is contained in the set of all integers with
$\Omega(n)\le q-3$, it follows from~\eqref{eq:LandauResult} that
the number of such integers $n \le x$ satisfies
\begin{equation}
\label{eq:N1} \cN_1(x)\ll
     \frac{x(\log_2x)^{q-4}}{\log x}\ .
\end{equation}

Next, let $n\in \cN_2$, and factor $n=p_1p_2\ldots p_k$, where
$p_1\le p_2\le \cdots \le p_{k}$ are primes, none of which is
equal to $q$; note that $k\ge q-2$. Let $a_j$ denote the residue
class of $p_j$ modulo $q$ for $j=1,\ldots, k$. For any nonempty
subset $S\subseteq \{1,\ldots,k\}$, $\prod_{j\in S} a_j$ is the
residue class of the divisor $d_S=\prod_{j\in S} p_j$ of~$n$.
Since $d_S\not\equiv 1\pmod q$, it follows that $k\le\kappa(G)$,
where $G$ is the abelian group $(\Z/q\Z)^*\cong\Z/(q-1)\Z$. Hence,
by Lemma~\ref{lem:comb1}, we have $k\le q-2$. Since $k\ge q-2$ for
each $n\in \cN_2$, it follows that $k=q-2$, and
Lemma~\ref{lem:comb1} further shows that $a_1\equiv \cdots\equiv
a_k\equiv a\pmod q$ for some primitive root $a$ modulo $q$.
Therefore, denoting by $U(q)$ the set of primitive roots modulo
$q$, we have
$$
\cN_2(x)=\sum_{a\in U(q)}{\cal Q}_{q-2}\big(x;q,\{a\}\big).
$$
Since $|U(q)|=\varphi(q-1)$, from Lemma~\ref{lem:anal1} we deduce
that
\begin{equation}
\label{eq:N2} \cN_2(x)=(1+o(1))\frac{\varphi(q-1)\,x(\log_2
x)^{q-3}}{(q-1)^{q-2}\,(q-3)!\,\log x}\ .
\end{equation}

Combining the estimates~\eqref{eq:N1} and~\eqref{eq:N2}, we have
\begin{equation}
\label{eq:Nstar} \cN^*(x)=(1+o(1))\frac{\varphi(q-1)\,x(\log_2
x)^{q-3}}{(q-1)^{q-2}\,(q-3)!\,\log x}\ .
\end{equation}
In view of the obvious relation
$$
\cN(x)=\sum_{\nu\ge 0}\cN^*(x/q^{\nu}),
$$
we see that
$$
\cN(x)=(1+o(1))(1-q^{-1})^{-1}\cN^*(x),
$$
which, together with~\eqref{eq:Nstar} yields the stated estimate
of Theorem~\ref{thm:easy}.

\subsection{Proof of Theorem~\ref{thm:hard}}

Fix the prime $q\ge 3$ and the integer $2\le a<q$, write $\cN$ for
$\cN(q,a)$, and let $\cN^*$ be the set of integers $n\in\cN$ that
are not divisible by $q$.

Throughout the proof, we fix a generator $g$ of the group
$G_r=\Z/p^r\Z$ with the property stated in
Lemma~\ref{lem:G/H(a)two}.  Here, $p^r=\paq$ as usual. We also
denote by $\phi_g:G_r\to G/H(a)$ and $\psi_g:G\to G_r$ the maps
defined in the statement and proof of Lemma~\ref{lem:G/H(a)two}.
Here, $G=(\Z/q\Z)^*$ as before.

For each $n\in\cN^*$, let $n=p_1\cdots p_k$ be a factorization of
$n$ as a product of primes, where $k=\Omega(n)$, and let
$\fM_n=\langle a_1,\ldots,a_k\rangle$ be the multiset in $G$ whose
elements are the congruence classes $p_j\pmod q$ for
$j=1,\ldots,k$.  As in the statement of Lemma~\ref{lem:technical},
we associate to $\fM_n$ a subgroup $H_n$ of $G$ and a multiset
$\fK_n\subseteq\fM_n$.

For every subgroup $H$ of $G$ with $a\not\in H$ and every multiset
$\fK$ in $G$, let $\cN_{H,\fK}$ denote the set of integers
$n\in\cN^*$, $n\le x$ such that $H_n=H$ and $\fK_n=\fK$. Our goal
is to estimate the number $\cN_{H,\fK}(x)$ of these, for every
pair $(H,\fK)$.

First, suppose that $H\ne H(a)$, and let $H$ and $\fK$ be fixed.
Put $y=\exp((\log x\log_3x)/\log_2x)$, and let
$$
\cN_1=\{n\in\cN_{H,\fK}~:~P(n)\le y\},
$$
where $P(n)$ denotes the largest prime factor of $n$. Using a
well--known result on {\it smooth numbers}; i.e., positive
integers $n$ whose largest prime factor is small with respect to
$n$ (see for example~\cite{CEP} or~\cite{Hilde}), we have
\begin{eqnarray}
\label{eq:N2T2} \cN_1(x)&\le& x\exp\(-(1+o(1))u\log u\)\nonumber\\
&=&\frac{x}{(\log x)^{1+o(1)}}=o\(\frac{x}{(\log x)^{1-1/\paq}}\),
\end{eqnarray}
where $u=(\log x)/\log y=(\log_2x)/\log_3x$.

Now let $\cN_2=\cN_{H,\fK}\setminus \cN_1$. For every integer $n$
in $\cN_2$, let $n=p_1\cdots p_k$ be a factorization of $n$ such
that $p_k=P(n)$, and put $m=p_1\cdots p_{k-1}$. For any fixed
value of $m$ obtained in this way, $p_k$ is a prime that satisfies
the inequalities
$$
x/m\ge p_k>y=\exp\(\frac{\log x\log_3x}{\log_2x}\);
$$
therefore, the number of possibilities for $p_k$ is at most
$$
\pi(x/m)\ll \frac{x}{m\log(x/m)}\le \frac{x\log_2 x} {m\log
x\log_3 x}.
$$
Note that $m=h_0k_0$, with
$$
h_0=\prod_{\substack{j=1\\a_j\in H}}^{k-1}p_j\mand
k_0=\prod_{\substack{j=1\\a_j\in\fK}}^{k-1}p_j,
$$
where each element $a_j\in G$ corresponds to the congruence class
$p_j\pmod q$ as before.  Then $h_0\in\cQ(q,H)$ in the notation of
Lemma~\ref{lem:A}, and we have $\Omega(k_0)\le |\fK|\le
L=(q-1)(q-3)$ by Lemma~\ref{lem:technical}$\,(i)$. Thus, summing
over the possible choices of $h_0$ and $k_0$, we see that
\begin{equation}
\label{eq:NHnotH(a)} \cN_2(x)\ll\frac{x\log_2 x}{\log x\log_3 x}
\Biggl(\sum_{\substack{h_0\in \cQ(q,H)\\h_0\le x}}
\frac{1}{h_0}\Biggr)\Biggl(\sum_{\substack{k_0\le
x\\\Omega(k_0)\le L}}\frac{1}{k_0}\Biggr).
\end{equation}
Using Lemma~\ref{lem:A} and partial summation, we derive the bound
\begin{equation}
\label{eq:**} \sum_{\substack{h_0\in \cQ(q,H)\\h_0\le
x}}\frac{1}{h_0}\ll
     (\log x)^{|H|/(q-1)}.
\end{equation}
On the other hand, we have
\begin{equation}
\label{eq:*}
\begin{split}
     \sum_{\substack{k_0\le x\\\Omega(k_0)\le L}}\frac{1}{k_0}
     &\ll\sum_{j\le L}\frac{1}{j!}
     \Biggl(\,\sum_{\substack{p\le x\\\nu\ge 1}}
     \frac{1}{p^{\nu}}\Biggr)^j\\
     &\ll\sum_{j\le L}\frac{1}{j!}\,(\log_2x+O(1))^j
     \ll (\log_2x)^L\ .
\end{split}
\end{equation}
Inserting the estimates~\eqref{eq:**} and~\eqref{eq:*}
into~\eqref{eq:NHnotH(a)}, it follows that
\begin{equation}
\label{eq:N3last} \cN_2(x)\ll\frac{x(\log_2x)^{L+1}}{(\log
x)^{1-|H|/(q-1)}\log_3 x}.
\end{equation}
Finally, by Lemma~\ref{lem:technical}$\,(ii)$, we have
$|H|<|H(a)|$ since $H\ne H(a)$ (and the group
$G=\left(\Z/q\Z\right)^*$ is cyclic). As $|H(a)|/(q-1)=1/\paq$ by
Lemma~\ref{lem:G/H(a)}, the estimates \eqref{eq:N2T2}
and~\eqref{eq:N3last} together imply that
\begin{equation}
\label{eq:NHnotH(a)2} \cN_{H,\fK}(x)=o\(\frac{x}{(\log
x)^{1-1/\paq}}\).
\end{equation}
Recall that the number of such pairs $(H,\fK)$ is bounded in terms
of $q$ so the above estimate is sufficient to easily absorb this
case into the error term.

It remains to consider the pairs with $H=H(a)$ and we turn our
attention to the problem of estimating $\cN_{H(a),\fK}(x)$ for a
fixed multiset $\fK$. In the case that $\fK=\varnothing$, it is
easy to see that
$$
\cN_{H(a),\varnothing}(x)=\cQ(x;q,H(a)).
$$
Hence, by Lemma~\ref{lem:A} we have
\begin{equation}
\label{eq:NH(a)Kempty}
\cN_{H(a),\varnothing}(x)=(1+o(1))\,\vartheta(q,H(a))
\frac{x}{(\log x)^{1-1/{\paq}}}.
\end{equation}

{From} now on, we assume that $\fK\ne\varnothing$.  We recall that
the inequality $|\fK|\le\paq-2$ holds by
Lemma~\ref{lem:technical}$\,(iii)$; in particular, $\paq\ge 3$ if
$\fK\ne\varnothing$.

First, suppose that $|\fK|<\paq-2$; note that this is possible
only if $\paq\ge 4$. For each $n\in\cN_{H(a),\fK}$, write
$n=h_0k_0$, where
$$
h_0=\prod_{\substack{j=1\\a_j\in H}}^{k}p_j\mand
k_0=\prod_{\substack{j=1\\a_j\in\fK}}^{k}p_j.
$$
Then $h_0\in\cS_1$ and $k_0\in\cS_2$, where
$$
\cS_1=\cQ(q,H(a))\mand
   \cS_2=\{n~:~\Omega(n)\le\paq-3\},
$$
and therefore,
$$
\cN_{H(a),\fK}(x)\le |\{(h_0,k_0)\in\cS_1\times\cS_2~:~h_0k_0\le
x\}|.
$$
Applying Lemma~\ref{lem:SSSSS} and making use of the estimates
provided by Lemma~\ref{lem:A} and~\eqref{eq:LandauResult}, we
obtain the bound
\begin{equation}
\label{eq:NH(a)Ksmall}
\cN_{H(a),\fK}(x)\ll\frac{x(\log_2x)^{\paq-3}}{(\log
x)^{1-1/\paq}},
\end{equation}
which again is of smaller order of magnitude than the main term
claimed by the theorem.

Now let $\fK$ be a multiset with cardinality $|\fK|=\paq-2$.
According to Lemma~\ref{lem:technical}$\,(iv)$, $\psi_g(\fK)$ is a
multiset $\fN_{p,r}(\eta,c)$ of the type considered in
Lemma~\ref{lem:listmultisets}; in other words,
$\fK\equiv\phi_g(\fN_{p,r}(\eta,c))\pmod{H(a)}$, or
$$
\fK=\big\langle g^ch_1,g^ch_2,\ldots,g^ch_\eta,
g^{\paq-c}h_{\eta+1},g^{\paq-c}h_{\eta+2},
\ldots,g^{\paq-c}h_{\paq-2}\big\rangle
$$
for some sequence $h_1,\ldots,h_{\paq-2}$ in $H(a)$.

For a fixed pair $(\eta,c)$, let $\cN_{\eta,c}$ be the disjoint
union
$$
\cN_{\eta,c}=\bigcup_{\psi_g(\fK)=\fN_{p,r}(\eta,c)}\cN_{H(a),\fK},
$$
and define the following subsets of $G$:
$$
G^+=\{g^ch:h\in H(a)\}\mand G^-=\{g^{\paq-c}h:h\in H(a)\}.
$$
For each $n\in\cN_{\eta,c}$, we can factor $n=h_0k_0l_0$, where
$$
h_0=\prod_{\substack{j=1\\a_j\in H(a)}}^{k}p_j,\qquad
k_0=\prod_{\substack{j=1\\a_j\in G^+}}^{k}p_j,\mand
l_0=\prod_{\substack{j=1\\a_j\in G^-}}^{k}p_j.
$$
Then $h_0\in\cS_1$, $k_0\in\cS_2$, and $l_0\in\cS_3$, where
\begin{eqnarray*}
\cS_1&=&\cQ(q,H(a)),\\
\cS_2&=&\cQ_\eta(q,G^+),\\
\cS_3&=&\cQ_\xi(q,G^-),
\end{eqnarray*}
with $\xi=\paq -2-\eta$. Conversely, if $h_0\in\cS_1$,
$k_0\in\cS_2$, and $l_0\in\cS_3$, and $n=h_0k_0l_0\le x$, then
$n\in\cN_{\eta,c}$.  Let us also define
$$
\cV=\{n~:~n=h_0l_0\text{~for
some~}h_0\in\cS_1\text{~and~}l_0\in\cS_3\}
$$
and
$$
\cW=\{n~:~n=h_0k_0l_0\text{~for
some~}h_0\in\cS_1,~k_0\in\cS_2,\text{~and~}l_0\in\cS_3\}.
$$
Then, since the sets $H(a)$, $G^+$ and $G^-$ are pairwise
disjoint, it is easy to see that the natural map
$\cS_1\times\cS_3\to\cV$ given by $(h_0,l_0)\mapsto h_0l_0$ is a
bijection.  Similarly, the natural map $\cV\times\cS_2\to\cW$
given by $(h_0l_0,k_0)\mapsto h_0k_0l_0$ is also a bijection. To
estimate $\cN_{\eta,c}(x)$, we apply Lemma~\ref{lem:SSSSS} twice:
first to the pair of sets $\cS_1$ and $\cS_3$, then to the pair of
sets $\cV$ and $\cS_2$.

By Lemma~\ref{lem:A}, we have
$$
\cS_1(x)=\cQ(x;q,H(a))=(1+o(1))\,\vartheta(q,H(a))
\,\frac{x}{(\log x)^{1-1/\paq}},
$$
and by Lemma~\ref{lem:anal1}, we have
$$
\cS_3(x)=\cQ_\xi(x;q,G^-)=(1+o(1))\frac{1}{\paq^\xi\,(\xi-1)!}
\frac{x(\log_2 x)^{\xi-1}}{\log x},
$$
where we have used the fact that $|G^-|= |H(a)|$. Applying
Lemma~\ref{lem:SSSSS} to the pair of sets $\cS_1$ and $\cS_3$, and
taking into account the bijection $\cS_1\times\cS_3\to\cV$
mentioned above, we get
\begin{equation*}
\begin{split}
\cV(x)&=|\{(h_0,l_0)\in\cS_1\times\cS_3~:~h_0l_0\le x\}|\\
&=(1+o(1))\frac{\vartheta(q,H(a))}{\paq^\xi\,\xi!}\frac{x(\log_2
x)^\xi}{(\log x)^{1-1/\paq}}.
\end{split}
\end{equation*}

To complete the estimate of $\cN_{\eta,c}(x)$, we must now
consider separately the cases $\eta=0$ and $\eta\ne 0$. Suppose
first that $\eta=0$ and $\xi=\paq-2$ (which can occur only if
$\paq$ is an odd prime; see Lemma~\ref{lem:listmultisets}). In
this case, $G^+=\varnothing$, $\cS_2=\{1\}$, and $\cW=\cV$;
consequently,
$$
\cN_{\eta,c}(x)=\cW(x)=(1+o(1))
\frac{\vartheta(q,H(a))}{\paq^{\paq-2}
\,(\paq-2)!}\frac{x(\log_2x)^{\paq-2}}{(\log x)^{1-1/\paq}}.
$$
Next, suppose that $\eta\ne 0$. By Lemma~\ref{lem:anal1}, we have
$$
\cS_2(x)=\cQ_\eta(x;q,G^+) =(1+o(1))\frac{1}{\paq^\eta\,(\eta-1)!}
\frac{x(\log_2 x)^{\eta-1}}{\log x}.
$$
Applying Lemma~\ref{lem:SSSSS} to the pair of sets $\cV$ and
$\cS_2$, and recalling the bijection $\cV\times\cS_2\to\cW$
described earlier, one has
\begin{equation*}
\begin{split}
\cN_{\eta,c}(x)=\cW(x)
&=|\{(h_0l_0,k_0)\in\cV\times\cS_2~:~h_0k_0l_0\le x\}|\\
&=(1+o(1))\frac{\vartheta(q,H(a))}{\paq^{\eta+\xi}\,\eta!\,\xi!}
\frac{x(\log_2x)^{\eta+\xi}}{(\log x)^{1-1/\paq}}.
\end{split}
\end{equation*}
Therefore, for all choices of $\eta$ and $c$, we obtain the
estimate
\begin{equation}
\label{eq:Netacmain} \cN_{\eta,c}(x)=(1+o(1))\binom{\paq-2}{\eta}
\frac{\vartheta(q,H(a))}{\paq^{\paq-2}\,(\paq-2)!}
\,\frac{x(\log_2x)^{\paq-2}}{(\log x)^{1-1/\paq}}.
\end{equation}

Taking into account the estimates \eqref{eq:N2T2},
\eqref{eq:NHnotH(a)2}, \eqref{eq:NH(a)Kempty},
\eqref{eq:NH(a)Ksmall} and \eqref{eq:Netacmain}, we find
$$
\cN^*(x)=\sum_{|\fK|=\paq-2}\cN_{H(a),\fK}(x)
+o\Biggl(\frac{x(\log_2x)^{\paq-2}}{(\log x)^{1-1/\paq}}\Biggr).
$$
Thus, if $\paq=2$, then $\cN_{H(a),\varnothing}(x)$ is the only
term in the above sum and
$$
\cN^*(x)=(1+o(1))\,\vartheta(q,H(a)) \frac{x}{(\log
x)^{1-1/{\paq}}}.
$$
If, on the other hand, $\paq\ge 3$, then
$$
\cN^*(x)=(1+o(1))\sum_{\eta,c} \binom{\paq-2}{\eta}
\frac{\vartheta(q,H(a))}{\paq^{\paq-2}\,(\paq-2)!}
\,\frac{x(\log_2x)^{\paq-2}}{(\log x)^{1-1/\paq}},
$$
where the sum runs over the possible values of $\eta$ and $c$
corresponding to the prime power $p^r=\paq$ (see
Lemma~\ref{lem:listmultisets}).  It is easy to see that
$$
\sum_{\eta,c} \binom{\paq-2}{\eta}=p^{r-1}\sum_{1\le j\le
p/2}\binom{p^r-2}{p^{r-1}j-1}
$$
holds for all possible values of $\paq$ (and both sides are equal
to $1$ if $\paq=2$); therefore, making use of
Lemma~\ref{lem:theta2V}, the definition~\eqref{eq:waq_defn}, and
the relation
$$
\cN(x)=(1+o(1))(1-q^{-1})^{-1}\cN^*(x),
$$
we obtain the estimate stated in the theorem.

\section{Concluding remarks}

We touch very briefly on a number of directions in which this work
might well be extended.

\smallskip
\noindent (1)\quad Further development of the main term in the
asymptotic formula: It is apparent that there are terms of only
slightly lower order in the asymptotic formula, some stepping down
by powers of $\log_2 x$ and others by powers of $\log x$. There
seems no reason why these could not be further elucidated although
a convenient description of the involved constants might be a lot
to expect.

\smallskip
\noindent (2)\quad Uniformity in the modulus: Certainly one can
trace through the above arguments to obtain results of this type.
If one wants however to obtain more than a very limited range of
applicability one would need to get at least some useful bounds
for the ``constants'' in the lower order main terms.

\smallskip
\noindent (3)\quad Subset avoidance: Rather than ask for the
number of integers whose divisors avoid a single residue class $a$
it seems natural to ask for the number of those whose divisors
avoid a subset $\cA$ of the reduced residue classes. Here, two
cases stand out as probably being quite similar to our existing
results, in the case that $\cA$ is a subgroup, to our first
theorem, and in the case that $\cA$ is a coset, to our second one.

\smallskip
\noindent (4)\quad General modulus: Although it could be combined
with any of the above, the removal of the restriction that the
modulus be prime is probably the most natural next step. In this
case it seems that little is needed beyond giving a count on the
number of different groups avoiding $a$ and having the same
maximal order, and then multiplying the previous result by this
number. It is clear that the contribution coming from integers
which correspond to more than one of these groups will give a
lower order of magnitude.  From the fundamental theorem for finite
abelian groups it is not hard to find a group-theoretic expression
for the number of such subgroups but to give this answer as an
explicit reasonable--looking function of the modulus may be a
different story.


\begin{thebibliography}{1000}

\bibitem{B} R. de la Bret\'eche, `R\'epartition des diviseurs dans les
progressions arithm\'etiques', \emph{Bull. \ London \ Math. \
Soc.} \textbf{32} (2000) no.~3, 257--262.

\bibitem{CEP}
E.~R.~Canfield, P.~Erd\H os and C.~Pomerance, `On a problem of
Oppenheim concerning ``factorisatio numerorum'',' \emph{J.\ Number
Theory} \textbf{17} (1983), no.~1, 1--28.

\bibitem{E} P. Erd\H os, `On the distribution of divisors of integers
in residue classes $\pmod d$', \emph{Bull.\ Math.\ Soc.\ Gr\'ece}
\textbf{6} Fasc. 1 (1965), 27--36.

\bibitem{E1} P. Erd\H os, `Some unconventional problems in number theory',
\emph{Ast\'erisque} \textbf{61} (1979), 73--82.

\bibitem{EH} P. Erd\H os and R.R. Hall, `Proof of a conjecture about the
distribution of divisors of integers in residue classes',
\emph{Math. \ Proc. \ Camb. \ Phil. \ Soc.} \textbf{79} (1976),
281--287.

\bibitem{ET} P. Erd\H os and G. Tenenbaum, `Ensemble de multiples
de suites finies',  Paul Erd\H os memorial collection.
\emph{Discrete \ Math.} \textbf{200}  (1999),  no.~ 1-3, 181--203.


\bibitem{Ha} R.R. Hall, `On some conjectures of Erd\H os in {\it Ast\'erisque}
I', \emph{J.\ Number\ Theory} \textbf{42}  (1992), 313--319.

\bibitem{Ha1} R.R. Hall, \emph{Sets of multiples}, Cambridge Texts in
Mathematics \textbf{118}, Cambridge U. Press, 1996.

\bibitem{HR} G.~H.~Hardy and S.~Ramanujan, `The normal number of prime
factors of an integer,' \emph{Quart.\ J.\ Math.\ (Oxford)}
\textbf{48} (1917), 76--92.

\bibitem{Hilde}
A.~Hildebrand, `On the number of positive integers $\leq x$ and
free of prime factors $>y$,'  \emph{J.\ Number Theory} \textbf{22}
(1986), no.~3, 289--307.

\bibitem{La}
E.~Landau, \emph{Handbuch der Lehre von der Verteilung der
Primzahlen.} Teubner, Leipzig-Berlin, 1909.

\bibitem{Nat}
M.~B.~Nathanson, \emph{Additive number theory. Inverse problems
and the geometry of sumsets.} Graduate Texts in Mathematics,
\textbf{165}. Springer-Verlag, New York, 1996.

\bibitem{Nat0}
M.~B.~Nathanson, \emph{Elementary methods in number theory.}
Graduate Texts in Mathematics, \textbf{195}. Springer-Verlag, New
York, 2000.

\bibitem{Uchiyama}
S.~Uchiyama, `On some products involving primes,' \emph{Proc.\
Amer.\ Math.\ Soc.} \textbf{28} (1971), 629--630.

\bibitem{Williams}
K.~S.~Williams, `Mertens' theorem for arithmetic progressions,'
\emph{J.\ Number Theory} \textbf{6} (1974), 353--359.

\bibitem{Wirs}
E.~Wirsing, `Das asymptotische Verhalten von Summen \"uber
multiplikative Funktionen,' (German) \emph{Math.\ Ann.}
\textbf{143} (1961), 75--102.

\end{thebibliography}
\end{document}